\input{style/arxiv-general.cfg}
\documentclass[MSNbibl,number,citesort,seceqn,rotating,dvips]{arxbj}
\makeatletter
   \@ifpackageloaded{graphicx}{}{\usepackage{graphicx}}
\makeatother
\usepackage{mathrsfs}
\usepackage{upgreek}
\usepackage{dcolumn}
\usepackage{graphicx}

\aid{0}
\volume{21}
\issue{2}
\pubyear{2015}
\firstpage{832}
\lastpage{850}
\doi{10.3150/13-BEJ589} 

\makeatletter
\newcommand{\rrvert}{\vert}
\newcommand{\llvert}{\vert}
\newcolumntype{d}[1]{D{.}{.}{#1}}
\renewcommand{\epsilon}{\varepsilon}
\newcommand{\ctan}{\operatorname{ctan}}
\newcommand{\arginf}{\operatorname{arg\,inf}}
\newcommand{\argsup}{\operatorname{arg\,sup}}
\newcommand{\fraca}[2]{{#1}/{#2}}
\newtheorem{them}{Theorem}[section]
\newproclaim{assumption}{Assumption}[section]
\newtheorem{corollary}{Corollary}[section]
\newtheorem{lemma}{Lemma}[section]
\makeatother

\begin{document}
\begin{frontmatter}

\title{Bayesian quantile regression with approximate likelihood}
\runtitle{Bayesian quantile regression with approximate likelihood}

\begin{aug}
\author[A]{\inits{Y.}\fnms{Yang} \snm{Feng}\thanksref{A}\ead[label=e1]{yfeng@google.com}},
\author[B]{\inits{Y.}\fnms{Yuguo} \snm{Chen}\corref{}\thanksref{B}\ead[label=e2]{yuguo@illinois.edu}} \and
\author[C]{\inits{X.}\fnms{Xuming} \snm{He}\thanksref{C}\ead[label=e3]{xmhe@umich.edu}}
\address[A]{Ads Metrics, Google Inc., Pittsburgh, PA 15206, USA. \printead{e1}}
\address[B]{Department of Statistics, University of Illinois at Urbana-Champaign, Champaign, IL 61820, USA. \printead{e2}}
\address[C]{Department of Statistics, University of Michigan, Ann Arbor, MI 48109, USA.\\ \printead{e3}}
\end{aug}

\received{\smonth{7} \syear{2012}}
\revised{\smonth{11} \syear{2012}}

%
\begin{abstract}
Quantile regression is
often used when a comprehensive relationship between a
response variable and one or more explanatory variables is desired. The
traditional frequentists' approach to quantile regression has been
well developed around asymptotic theories and efficient
algorithms. However, not much work has been published under the
Bayesian framework. One challenging problem for Bayesian
quantile regression is that the full likelihood has no parametric forms.
In this paper, we propose a Bayesian quantile
regression method, the linearly interpolated density (LID) method,
which uses a linear interpolation of the quantiles to approximate the
likelihood.
Unlike most of the existing methods that aim at tackling one
quantile at a time, our proposed method estimates the
joint posterior distribution of multiple quantiles, leading to
higher global efficiency for all quantiles of interest.
Markov chain Monte Carlo algorithms are developed to carry out
the proposed method.
We provide convergence results that justify both the algorithmic convergence
and statistical approximations to an integrated-likelihood-based posterior.
From the simulation results, we verify that LID
has a clear advantage over other existing methods in estimating
quantities that relate to two or more quantiles.
\end{abstract}

%

\end{frontmatter}

\section{Introduction}
\label{intro}

Quantile regression, as a supplement to the mean regression, is often
used when a comprehensive relationship between the response variable
$y$ and the explanatory variables $x$ is desired. Consider the
following linear model:
%
%
\begin{equation}
\label{model1} y_i=x_i^{T}\beta+
\epsilon_i,\qquad i=1,2,\ldots,n,
\end{equation}
where $y_i$ is the response variable, $x_i$ is a $p\times1$ vector
consisting of $p$ explanatory variables, $\beta$ is a $p\times1$
vector of coefficients for the explanatory variables, and
$\epsilon_i$ is the error term. The quantile regression analysis models
the $\tau$th conditional quantile of $y$ given $x$ as:
%
%
\begin{equation}
\label{model2} Q_{y_i}(\tau|x_i)=x_i^{T}
\beta(\tau),\qquad i=1,2,\ldots,n,
\end{equation}
which is equivalent to (\ref{model1}) with $Q_{\epsilon_i}(\tau|x_i)=0$.
The $\tau$-specific coefficient vector
$\beta(\tau)$ can be estimated by minimizing the loss function:
%
%
\begin{equation}
\label{QR_objective_function} \min_{\beta(\tau)}\sum_{i=1}^n
\rho_{\tau}\bigl(y_i-x_i^{T}\beta(
\tau)\bigr),
\end{equation}
where $\rho_{\tau}(u)=u\tau$ if $u\geq0$, and $\rho_{\tau
}(u)=u(\tau -1)$ if $u<0$;
see Koenker \cite{Koe05}.

To make inference on the quantile regression, one could
use the asymptotic normal distribution of the estimates or use the
bootstrap method. Aside from the regular bootstrap such as the residual
bootstrap and the
$(x,y)$ bootstrap, one could also use
Parzen, Wei and Ying \cite{Par94}'s method or the
Markov chain marginal bootstrap method (He and Hu \cite{he02}).

In contrast to the rich literature on quantile regression with the frequentist view,
not much work has been done under the Bayesian framework.
The most challenging problem for Bayesian quantile
regression is that the likelihood is usually not available unless
the conditional distribution for the error is assumed.

Yu and Moyeed \cite{Yu01} proposed an idea of employing a
likelihood function based on the asymmetric Laplace distribution. In
their work, Yu and Moyeed assumed that the error term follows an
independent asymmetric Laplace distribution
%
%
\begin{equation}
f_{\tau}(u)=\tau(1-\tau)\mathrm{e}^{-\rho_{\tau}(u)},\qquad u\in{R},
\end{equation}
where $\rho_{\tau}(u)$ is the loss function of quantile regression.
The asymmetric Laplace distribution is very closely related to
quantile regression since the mode of $f_{\tau}(u)$ is the solution
to (\ref{QR_objective_function}).
Reich, Bondell and Wang \cite{Reich10} developed a Bayesian
approach for quantile regression
assuming that the error term follows an infinite mixture of Gaussian
densities and their prior for the residual density is stochastically
centered on the asymmetric Laplace distribution.
Kottas and Gelfand \cite{Kot01} implemented a Bayesian median
regression by
introducing two families of distributions with median zero and the
Dirichlet process prior.
Dunson and Taylor \cite{Dun05} used a substitution likelihood
proposed by
Lavine \cite{Lav95} to make inferences based on the
posterior distribution.
One property of
Dunson and Taylor's method is that it allows regression on multiple
quantiles simultaneously.
Tokdar and Kadane \cite{Tokdar12} proposed a semiparametric
Bayesian approach for
simultaneous analysis of quantile regression models based on the
observation that when there is only a univariate covariate, the
monotonicity constraint can be satisfied by interpolating two monotone
curves, and the Bayesian inference can be carried out by specifying a
prior on the two monotone curves.
Taddy and Kottas \cite{Tad09}
developed a fully
nonparametric model-based quantile regression based on Dirichlet
process mixing. Kottas and Krnjaji{\'c} \cite{Kottas09} extended
this idea to the case where
the error distribution changes
nonparametrically with the covariates.
Recently, Yang and He \cite{Yang11} proposed a Bayesian
empirical likelihood method which targets on estimating multiple
quantiles simultaneously, and justified the validity of the posterior
based inference.

In this paper, we propose a Bayesian method, which aims at estimating
the joint posterior distribution of multiple quantiles and achieving
``global'' efficiency for quantiles of interest.
We consider a Bayesian approach to estimating multiple quantiles as
follows. Let $\tau_1,\ldots,\tau_m$ be $m$ quantiles in model (\ref
{model2}) and $B_m=(\beta(\tau_1),\ldots,\beta(\tau_m))$.
Let $X=(x_1,\ldots,x_n)'$ and $Y=(y_1,\ldots,y_n)$ be the observations
of size $n$.
For each pair of observation $(x_i,y_i)$, the likelihood $L(B_m|x_i,
y_i)=p(y_i|x_i,B_m)$
is not available. However if we include $f_i$, the probability density
function (pdf) of the conditional distribution $y|x_i$, as the nuisance
parameter, then the likelihood $L(B_m, f_i|x_i, y_i)=p(y_i|x_i,B_m,
f_i)=f_i(y_i)$. This is to treat Bayesian quantile regression as a
semi-parametric problem: the parameter of interest is finite
dimensional and the nuisance parameter is nonparametric. To eliminate
the nuisance parameter, we use the integrated likelihood methods
recommended by Berger, Liseo and Wolpert \cite{Berger99}. More specifically,
let $\theta_{f_i}$ be all the quantiles of $f_i$, and
$\theta_{m,i}=x_iB_m$ be the $m$ quantiles of interest. We can define
$p(y_i|x_i,B_m)$ as
%
%
\begin{equation}
\label{def} p(y_i|x_i, B_m)=\int
_{f_i\in\mathscr{F}_{\theta_{m,i}}} p(y_i|\theta_{f_i})\,\mathrm{d}
\Pi_{\theta_{m,i}}(f_i),
\end{equation}
where
$\mathscr{F}_{\theta_{m,i}}$ denotes the subset of
well-behaved pdfs (will be defined precisely in Section~\ref{sec3.2})
with those $m$ quantiles equal to $\theta_{m,i}$, $\Pi_{\theta
_{m,i}}(\cdot)$ denotes the prior on $f_i|\theta_{m,i}\in\mathscr
{F}_{\theta_{m,i}}$ (will be specified in Section~\ref{sec3.2}), and
$p(y_i|\theta_{f_i})=f_i(y_i)$ because
$f_i (y |x_i)$ is determined by the conditional quantile functions.
Here, $p(y_i|x_i, B_m)$ can be viewed as an integral of a function or
an expectation with the densities as the random variable.
The posterior distribution of $B_m|X,Y$ can be written as
%
%
\begin{equation}
\label{target distribution} p(B_m|X,Y)\propto\pi_m(B_m|X)L(Y|X,B_m),
\end{equation}
where $\pi_m(B_m|X)$ is the prior on $B_m$ and $L(Y|X,B_m)=\prod_{i=1}^n p(y_i|x_i, B_m)$.

One practical difficulty with the above approach is that the
integration step to remove the nuisance parameter is computationally
infeasible except for the case of $m=1$ (Doss \cite{Doss85}). To circumvent
this issue, we consider a different approximation to the likelihood.
Note that $x_iB_m$ gives the $m$ quantiles of the conditional
distribution $y|x_i$ based on model (\ref{model2}). These $m$ quantiles
can be used to construct an approximate conditional distribution
$y|x_i$ through linear interpolation.
With this approximate likelihood, an approximate posterior distribution
becomes available. We show that the total variation distance between
the approximate posterior distribution and $p(B_m|X,Y)$ (the posterior
based on the integrated likelihood) goes to 0 as $\tau_1,\ldots,\tau_m$
becomes dense in $(0,1)$ as $m\rightarrow\infty$. A Markov chain Monte
Carlo (MCMC) algorithm can then be developed to sample from the
approximate posterior distribution.
The recent work of Reich, Fuentes and Dunson \cite{Rei11} used
large-sample approximations to
the likelihood to do Bayesian quantile regression. Their approach also
aims to achieve global efficiency over multiple quantiles, and can
adapt to account for spatial correlation. In contrast, our work uses
approximations at a fixed sample size $n$ and provides a Bayesian
interpretation of the posterior quantities.

The rest of the paper
is organized as follows. Section~\ref{sec2} introduces the proposed
method. Section~\ref{sec3} provides the convergence property of the
algorithm as well as the approximate posterior distribution.
Section~\ref{sec4} compares the proposed method with some existing methods
through simulation studies and applies the proposed method to real
data. Section~\ref{sec6} provides concluding remarks.

\section{Methodology}
\label{sec2}

In this section, we describe the linearly interpolated density to be
used in approximating the likelihood, and then give the layout of our
MCMC algorithm for posterior inference.
We list again the basic setting introduced in Section~\ref{intro}.
Let $X=(x_1,\ldots,x_n)'$ and $Y=(y_1,\ldots,y_n)$ be the observations.
Let $\tau_1,\ldots,\tau_m$ be $m$ quantiles in model (\ref{model2}) and
$B_m=(\beta(\tau_1),\ldots,\beta(\tau_m))$.
We are interested in the posterior distribution $B_m|X,Y$.

\subsection{Linearly interpolated density}

The likelihood is generally not assumed under the
quantile regression model, but $x_iB_m$ gives the $m$ quantiles of the
conditional distribution $y|x_i$.
With the linearly interpolated density based on the $m$ quantiles, we
can approximate the true likelihood from a sequence of specified
quantile functions.

Here is how the linear interpolation idea works in a simple setting.
Suppose $Z\sim F(z)$, where $F(z)$ is the cumulative distribution
function (cdf) of $Z$. Let $f(z)$ be
the pdf of $Z$. Let $\tau_z=F(z)$, and
$\tau_1,\tau_2$ be two constants such that
$0\le\tau_1<\tau_z<\tau_2\le1$. Then
$F^{-1}(\tau_1)<z<F^{-1}(\tau_2)$ if $f(z)$ is continuous
and non-zero on the support of $Z$. We can approximate $f(z)$ by
%
%
\begin{equation}
\frac{\tau_2-\tau_1}{F^{-1}(\tau_2)-F^{-1}(\tau_1)},
\end{equation}
because
%
%
\begin{equation}
\frac{\tau_2-\tau_1}{F^{-1}(\tau_2)-F^{-1}(\tau_1)}=\frac{\tau
_2-\tau
_1}{\frac{\mathrm{d}}{\mathrm{d}\tau}F^{-1}(\tau^{*})(\tau_2-\tau_1)}=f\bigl(z^*\bigr),
\end{equation}
where $\tau_1<\tau^*<\tau_2$ and
$z^*=F^{-1}(\tau^*)\in(F^{-1}(\tau_1),F^{-1}(\tau_2))$.

Now we extend the interpolation idea to model (\ref{model2}). Given
$B_m=(\beta(\tau_1), \beta(\tau_2),\ldots,\break \beta(\tau_m))$, we could
calculate the
linearly interpolated density $\hat{f}_i(y_i|x_i,B_m)$, $i=1,2,\ldots
,n$, by
%
%
\begin{eqnarray}
\label{den} \hat{f}_{i}(y_i|x_i,B_m)&=&
\Biggl[\sum_{j=1}^{m-1}I_{\{y_i\in
(x_i\beta
(\tau_j),  x_i\beta(\tau_{j+1}))\}}
\frac{\tau
_{j+1}-\tau
_j}{x_i\beta(\tau_{j+1})-x_i\beta(\tau_j)} \Biggr]
\nonumber
\\[-8pt]
\\[-8pt]
&&{}+I_{\{y_i\in(-\infty, x_i\beta(\tau_1))\}}\tau_1f_1(y_i)
+I_{\{y_i\in(x_i\beta(\tau_m), \infty)\}}(1-\tau_m)f_2(y_i),
\nonumber
\end{eqnarray}
where $f_1$ is distributed as the left half of $N(x_i\beta(\tau
_1),\sigma^2)$, $f_2$ is distributed as the right half of $N(x_i\beta
(\tau_m),\sigma^2)$, and $\sigma^2$ is some pre-specified parameter.

Let $\hat{p}_m(Y|X,B_m)=\prod_{i=1}^n\hat{f}_{i}(y_i|x_i,B_m)$
denote the
approximate likelihood. One possible prior $\pi_m(B_m|X)$ on $B_m$ is a
truncated normal $N(\mu,\Sigma)$ satisfying
%
%
\begin{equation}
\label{cons1} x_i\beta(\tau_1)<x_i\beta(
\tau_2)<\cdots<x_i\beta(\tau_m),\qquad
i=1,2,\ldots,n.
\end{equation}
Since we include the intercept in model (\ref{model2}),
the first element of $x_i$ is 1, and at least the parallel quantile
regression lines satisfy (\ref{cons1}).
The corresponding posterior is
%
%
\begin{equation}
\label{approx-posterior} \hat{p}_m(B_m|X,Y)= \frac{\pi_m(B_m|X)\hat{p}_m(Y|X,B_m)}{\hat{p}_m(Y|X)},
\end{equation}
where $\hat{p}_m(Y|X)=\int\pi_m(B_m|X)\hat{p}_m(Y|X,B_m)\,\mathrm{d}B_m$.
In the next section, we give a MCMC algorithm to sample $B_m$ from this
posterior.
We show later that the total variation distance between this posterior
distribution and the target posterior $p(B_m|X,Y)$ goes to 0 as $m$
goes to infinity.

\subsection{Algorithm of the linearly interpolated density (LID) method}
\label{algorithm}

We incorporate the linearly interpolated density into the following
modified Metropolis--Hastings algorithm to draw samples from $\hat
{p}_m(B_m|X,Y)$.
\begin{enumerate}[6.]
\item[1.] Choose an initial value $B_m^0$ for $B_m$. One good choice is
to use the parallel quantile estimates, that is, all the slopes for
the quantiles are the same and the intercepts are different. We could
use the \emph{quantreg} (a function in R) estimates of the slopes for
the median as the initial slopes, and use
the \emph{quantreg} estimates of the intercepts for each quantile as
the initial intercepts. In case a lower quantile has a larger intercept
than an upper quantile, we could order the intercepts such that the
intercepts increase with respect to $\tau$. If there are ties, we could
add an increasing sequence with respect to $\tau$ to the intercepts to
distinguish them. Another possible choice for the initial value is to
use Bondell, Reich and Wang \cite{Bond10}'s estimate which
guarantees the non-crossing of the quantiles.

\item[2.] Approximate the densities. With the initial values of the
parameters, we can calculate the linearly interpolated density $\hat
{f}^0_{i}(y_i|x_i, B_m^0)$, $i=1,2,\ldots,n$, by
plugging $B^0_m$ into equation (\ref{den}).
Let
$L^{0}=\prod_{i=1}^n \hat{f}^0_{i}(y_i|x_i, B_m^0)$.

\item[3.] Propose a move. Suppose we are at the $k$th iteration. Randomly
pick a number $\tau_j$ from $\tau_1,\tau_2,\ldots,\tau_m$ and then
randomly pick a component $\beta_l^{k-1}(\tau_j)$ of $\beta
^{k-1}(\tau_j)$
to update. To make sure that the proposed point $\beta_l^*(\tau_j)$
satisfies constraint (\ref{cons1}), we can calculate a lower bound
$l_{j,l}$ and an upper bound $u_{j,l}$ for $\beta_l^*(\tau_j)$ and
generate a
value for $\beta_l^*(\tau_j)$ from $\operatorname{Uniform}(l_{j,l},u_{j,l})$. In case
$l_{j,l}=-\infty$
or $u_{j,l}=\infty$, we will use
a truncated normal as the proposal distribution. The details on how to
find the bounds are in Appendix~\ref{bound}.
Denote $\beta^{*}(\tau_j)$ as the updated $\beta^{k-1}(\tau_j)$ by
replacing its $l$th component $\beta_l^{k-1}(\tau_j)$ by the proposed
value $\beta_l^*(\tau_j)$.

\item[4.] Set
$B^*_m=(\beta^{k-1}(\tau_1),\ldots,\beta^{k-1}(\tau_{j-1}),\beta
^{*}(\tau_j),\beta^{k-1}(\tau_{j+1}),\ldots,\beta^{k-1}(\tau_m))$.
We can calculate the linearly interpolated density
$\hat{f}^*_i(y_i|x_i, B_m^*)$, $i=1,2,\ldots,n$, by
plugging $B^*_m$ into equation (\ref{den}).
Let $L^*=\prod_{i=1}^n \hat{f}^*_i(y_i|x_i, B_m^*)$.

\item[5.] Calculate the acceptance probability
%
%
\begin{equation}
r=\min \biggl(1,\frac{\pi_m(B^*_m|X)L^*q(B^{*}_m\to
B^{k-1}_m)}{\pi_m(B^{k-1}_m|X)L^{k-1}q(B^{k-1}_m\to B^{*}_m)} \biggr),
\end{equation}
where $q(B^{k-1}_m\to B^*_m)$ denotes the transition probability
from $B^{k-1}_m$ to $B^*_m$. Notice that
these two transition probabilities cancel out if we choose
symmetric proposals. Let $B^k_m=B^*_m$ with probability $r$, and
$B^k_m=B^{k-1}_m$ with probability $1-r$. If
$B^k_m=B^*_m$, then $L^k=L^*$; otherwise $L^k=L^{k-1}$.

\item[6.] Repeat steps 3--5 until the desired number of iterations is
reached.
\end{enumerate}

\section{Theoretical properties}
\label{sec3}

In this section, we give the stationary distribution of the Markov
chain in Section~\ref{algorithm} for fixed $m$, and study the limiting
behavior of the stationary distribution
as $m\rightarrow\infty$.

\subsection{Stationary distribution}

Since we replace the true probability density function by the linearly
interpolated density in the Metropolis--Hastings algorithm in
Section~\ref{algorithm}, it is not obvious what the stationary
distribution of
the Markov chain is. The following theorem, whose proof is in Appendix~\ref{proof1}, says that the Markov chain converges to $\hat
{p}_m(B_m|X,Y)$ defined in (\ref{approx-posterior}).

%
\begin{them}
\label{theo}
The stationary distribution of the Markov chain constructed in
Section~\ref{algorithm} is
$\hat{p}_m(B_m|X,Y)$.
\end{them}

This theorem implies that we can use the algorithm in Section~\ref
{algorithm} to draw samples from $\hat{p}_m(B_m|X,Y)$.

\subsection{Limiting distribution}
\label{sec3.2}

In this section, we show that as $m\to\infty$, the total variation
distance between the
stationary distribution $\hat{p}_m(B_m|X,Y)$
and the target distribution $p(B_m|X,Y)$ (defined in (\ref{target
distribution}))
goes to 0.
The proof requires the following assumption about $f_i$, the
probability density function of the conditional distribution $y|x_i$.
All the results are stated for a given sample size $n$.

\begin{assumption}
\label{assump1}
Let $q_{f,\tau}$ be the $\tau$th quantile of $f$, and $M_1$,
$M_2$ and $c$ be constants. The densities of $y|x_i$ are in the set
$\mathscr{F}=\{f|\int f \,\mathrm{d}x=1,0\le f\leq M_1, |f'|<M_2, \mbox{and }
f(x)<c/\sqrt{m}\mbox{ for }x<q_{f,1/m}\mbox{ and for }
x>q_{f,(m-1)/m}, m=2,3,\ldots\}$.
\end{assumption}

The assumption implies that $\mathscr{F}$ is a
set of bounded probability density functions with bounded
first derivatives and controlled tails. The restrictions on the tails
are not hard to satisfy. The
Cauchy distribution, for example, is in the set.
For the Cauchy distribution, the $\frac{1}{m}$th quantile is
$q_{\fraca{1}{m}}=\tan(\uppi(\frac{1}{m}-\frac{1}{2}))=-\ctan(\frac{\uppi}{m})$,
so
$f(q_{\fraca{1}{m}})=\frac{1}{\uppi}\frac{1}{1+\ctan^2(\fraca{\uppi
}{m})}=\frac
{1}{\uppi}\sin^2(\frac{\uppi}{m})=\mathrm{O}(\frac{1}{m^2})<\frac{c}{\sqrt{m}}$
for some $c$.
The set $\mathscr{F}_{\theta_{m,i}}$ appeared in (\ref{def}) denotes
the subset of $\mathscr{F}$ that
contains all the pdfs with those $m$ quantiles equal to $\theta_{m,i}=x_iB_m$.

We now specify the prior on $f_i(\cdot|x_i)\in\mathscr{F}$, denoted by
$\Pi(f_i)$, and the prior on $f_i|\theta_{m,i}\in\mathscr
{F}_{\theta
_{m,i}}$, denoted by $\Pi_{\theta_{m,i}}(f_i)$. We know from
(\ref{model2}) that the $\tau$th quantile of $f_i(\cdot|x_i)$, the
conditional distribution of $y$ given $x=x_i$, is $x_i^{T}\beta(\tau)$.
Let us consider $\beta(\tau)$ as a function of $\tau$, where
$0\leq\tau\leq1$. Because $x_i^{T}\beta(\tau)$, $0\leq\tau\leq
1$, determines
all the quantiles of $f_i(\cdot|x_i)$ based on (\ref{model2}), and
therefore determines $f_i(\cdot|x_i)$ (Koenker \cite{Koe05}), the prior on
$f_i(\cdot|x_i)$ can be induced
from the prior on $\beta(\tau)$.
To satisfy Assumption~\ref{assump1}, we use a Gaussian process prior on
$\beta''(\tau)$ so
that $\beta(\tau)$ has the second derivative, and then $f_i$'s have the
first derivative.
The prior $\Pi(f_i)$ on $f_i(\cdot|x_i)$ is induced from the prior on
$\beta(\tau)$. The prior $\Pi_{\theta_{m,i}}(f_i)$ on $f_i|\theta
_{m,i}$ is
induced by $\Pi(f_i)$.
The prior on $B_m$ can be obtained from the
prior on $\beta(\tau)$, because $B_m$ is a vector of $m$ points on
$\beta(\tau)$.
With the specification of these priors, $p(y_i|x_i, B_m)$ and
$p(B_m|X,Y)$ given in
(\ref{def}) and (\ref{target distribution}) are well-defined.

To study the limiting distribution as $m\to\infty$, we assume the
sequence of quantile levels satisfies the following condition:
%
%
\begin{equation}
\label{condition} \Delta\tau=\max_{0\le j\le
m}(\tau_{j+1}-
\tau_j)=\mathrm{O}\biggl(\frac{1}{m}\biggr),
\end{equation}
where $\tau_0=0$ and $\tau_{m+1}=1$.
This condition is not difficult to satisfy. For example, we can start from
$m_0=M_0$ quantile levels: $\tau=\frac{1}{M_0+1}, \frac{2}{M_0+1},
\ldots, \frac{M_0}{M_0+1}$, which include the quantiles of interest.
We add new $\tau$'s one
by one so that the new $\tau$ divides one of the previous intervals in
halves, that is, $\tau=\frac{1}{2(M_0+1)},
\frac{3}{2(M_0+1)},\ldots, \frac{2M_0+1}{2(M_0+1)},
\frac{1}{4(M_0+1)}, \frac{3}{4(M_0+1)},\ldots,\break
\frac{4M_0+3}{4(M_0+1)}$ and so on. For this sequence of quantiles,
we have
$\Delta\tau=\max_{0\le j\le m}(\tau_{j+1}-\tau_j)\leq\frac
{2}{m}=\mathrm{O}(\frac
{1}{m})$.

To prove the convergence of distributions, we
use the total variation norm,
$\Vert \mu_1-\mu_2\Vert _{\mathrm{TV}}=\sup_A|\mu_1(A)-\mu_2(A)|$ for two probability
measures $\mu_1$ and $\mu_2$, where $A$ denotes any measurable set.
It is more convenient to use the following equivalent definition
(Robert and Casella \cite{Rob04}, page~253): $\Vert \mu_1-\mu
_2\Vert _{\mathrm{TV}}=\frac{1}{2}\sup_{|h|\leq1}|\int
h(x)\mu_1(\mathrm{d}x)-\int h(x)\mu_2(\mathrm{d}x)|$.
The following theorem
gives the limiting distribution of the stationary distribution as
$m\to\infty$.

%
\begin{them}
\label{theo1}
$\|\hat{p}_m(B_m|X,Y)-p(B_m|X,Y)\|_{\mathrm{TV}}\rightarrow0$ as
$m\rightarrow\infty$, assuming $\tau_{j+1}-\tau_j=\mathrm{O}(\frac{1}{m})$.
\end{them}

The proof is in Appendix~\ref{proof2}.
As a consequence of Theorem~\ref{theo1}, we have
the following corollary.

\begin{corollary}
Let $\eta$ be the quantiles of interest, which is contained in $B_m$.
We have $\|\hat{p}_m(\eta|X,Y)-p(\eta|X,Y)\|_{\mathrm{TV}}\to0$ as $m\to
\infty$,
assuming $\tau_{j+1}-\tau_j=\mathrm{O}(\frac{1}{m})$.
\end{corollary}

The above corollary says that by the linearly interpolated density
approximation the posterior distribution of the quantiles of interest
converges to the target distribution. The theorem requires that we need
to increase $m$ in the algorithm. Although $m$ is fixed in
applications, the convergence   result lends support to $\hat
{p}_m(B_m|X,Y)$ as an approximation.

\section{Comparison of LID with other methods}
\label{sec4}

In this section, we compare the proposed method with some existing
methods through three simulation studies.
In the quantile regression model (\ref{model2}), if the conditional
densities $f_i(y|x_i)$ are different for different observation $i$, one could
apply weighted quantile regression to improve the efficiency of
estimates (Koenker \cite{Koe05}, page 160). In this
case, the loss function would be:
%
%
\begin{equation}
\min_{\beta(\tau)}\sum_{i=1}^nw_i
\rho_{\tau}\bigl(y_i-x_i^{T}\beta (
\tau)\bigr),
\end{equation}
where $w_i$ denotes the weight for the $i$th observation. The optimal
weight is the conditional density $f_i(y|x_i)$ at the $\tau$th
quantile. Because
the density is not available generally, one could approximate the
density by a nonparametric density estimate. One simple way
is to use
%
%
\begin{equation}
\hat{w}_i=\frac{2\Delta\tau}{x_i^{T}(\beta^{rq}(\tau+\Delta\tau
)-\beta
^{rq}(\tau-\Delta\tau))},\qquad i=1,2,\ldots,n,
\end{equation}
where $\beta^{rq}$ denotes the unweighted quantile regression estimate.
When the weight is negative due to crossing of quantile estimates, we
just set the weight to be 0. This occurs with probability tending to 0
as $n$ increases.
To make inference, one could
use the asymptotic normal distribution of the estimates or use the
bootstrap method.

\subsection{Example~1}

The data were generated from the following model
%
%
\begin{equation}
\label{simulation1} y_i=a+bx_i+(1+x_i)
\epsilon_i,\qquad i=1,2,\ldots,n,
\end{equation}
where $\epsilon_i$'s are independent and identically distributed
(i.i.d.) as $N(0,1)$. We chose $n=100$, \mbox{$a=5$} and $b=1$. The covariate
$x_{i}$ was generated from $\operatorname{lognormal}(0,1)$.
The corresponding quantiles of interest are
%
%
\begin{equation}
Q_{y_i}(\tau|x_i)=a(\tau)+b(\tau)x_i,\qquad
i=1,2,\ldots,n, \tau=\frac
{1}{m+1},\ldots,\frac{m}{m+1}.
\end{equation}
Here we report the results on the $0.25$, $0.5$ and $0.75$ quantiles
and the difference between the $0.75$ and $0.5$ quantiles by comparing
the mean squared error (MSE) for the slope estimates from five
different methods: the proposed linearly interpolated density method
(LID), the regular regression of quantiles (RQ), the weighted RQ with
estimated weights (EWRQ) (Koenker \cite{Koe05}),
the pseudo-Bayesian method of Yu and Moyeed \cite{Yu01},
and the approximate Bayesian
method of Reich, Fuentes and Dunson \cite{Rei11}.
We generated 100 data sets for computing the MSE.

For LID and Yu and Moyeed's method, we used the normal prior $N(0,100)$
for each parameter $a(\tau)$ and $b(\tau)$. For LID, we chose $m=49$,
equally spaced quantiles between 0 and 1 (which include the quantiles
of interest: $0.25$, $0.5$ and $0.75$), and the length of the Markov
chain is 1\,000\,000 (half of the samples were used as burn-in). We ran
such a long chain because we updated 98 parameters one at a time, which
means we updated each parameter about 10\,000 times on average. Every
thousandth sample in the chain is taken for the posterior inference.
For Yu and Moyeed's method, a Markov chain with length 5\,000 (half of
the samples were used as burn-in) seems enough for the inference,
partially because Yu and Moyeed's method is dealing with one quantile
at a time and has only two parameters. For Reich et al.'s method, we
simply used their code and set the length of the chain to be 2\,000
(half of the samples were used as burn-in).
Notice that for LID and Reich et al.'s method, only one run is needed
to provide all results in the table, and other methods have to run for
each $\tau$.

From the results in Table~\ref{t100}, we can see that LID did better
than RQ and Yu and Moyeed's method. Comparing with weighted RQ and
Reich et al.'s method, LID gave better estimates for upper quantiles
but poorer estimates for lower quantiles.
For estimating the differences of quantiles, LID is clearly the best
among all the methods.

\begin{table}
\tablewidth=\textwidth
\tabcolsep=0pt
\caption{$n\times\mathrm{MSE}$ and its standard error (in parentheses) for
Example~1}\label{t100}
\begin{tabular*}{\textwidth}{@{\extracolsep{\fill}}lllll@{}}
\hline
Methods & $b(0.25)$ & $b(0.5)$ & $b(0.75)$ & $b(0.75)-b(0.5)$
\\
\hline
RQ & 23 (4) & 19 (2) & 19 (3) & 15 (2)
\\
EWRQ & 16 (2) & 13 (2) & 15 (3) & 11 (2)
\\
LID & 22 (4) & 15 (2) & 13 (1) & \hphantom{0}3 (0.6)
\\
Yu and Moyeed & 21 (4) & 17 (2) & 16 (3) & 10 (1)
\\
Reich et al. & 16 (2) & 15 (2) & 23 (3) & 11 (1)
\\
\hline
\end{tabular*}
\end{table}

\subsection{Example~2}

The data were generated from the following model
%
%
\begin{equation}
\label{simulation2} y_i=a+bx_{1,i}+cx_{2,i}+(1+x_{1,i}+x_{2,i})
\epsilon_i,\qquad i=1,2,\ldots,n,
\end{equation}
where $\epsilon_i$'s are i.i.d. from $N(0,1)$. In the simulations, we
chose $n=100$, $a=5$, $b=1$, and $c=1$. The covariates $x_{1,i}$
was generated from $\operatorname{lognormal}(0,1)$ and $x_{2,i}$ was generated from
$\operatorname{Bernoulli}(0.5)$.
The corresponding quantiles of interest are
%
%
\begin{equation}
Q_{y_i}(\tau|x_i)=a(\tau)+b(\tau)x_{1,i}+c(
\tau)x_{2,i},\qquad i=1,2,\ldots ,n, \tau=\frac{1}{m+1},\ldots,
\frac{m}{m+1}.
\end{equation}

We compared the five methods with the same performance criterion as
Example~1.
We generated 400 data sets for computing the MSE.
The results are in Table~\ref{t13}. We see that for the quantile
estimates, LID (with $m=15$) and EWRQ perform similarly, and LID
outperforms RQ and Yu and Moyeed's method.
For estimating the difference between quantiles, LID outperforms RQ,
EWRQ, and Yu and Moyeed's method.
Comparing with Reich et al.'s method, LID gave better estimates for
parameter $b$ but poorer estimates for parameter $c$.

\begin{table}[t]
\tablewidth=\textwidth
\tabcolsep=0pt
\caption{$n\times\mathrm{MSE}$ and its standard error (in parenthesis) for
Example~2}\label{t13}
\begin{tabular*}{\textwidth}{@{\extracolsep{\fill}}lllllll@{}}
\hline
Methods & $b(0.5)$ & $b(0.75)$ & $b(0.75)-b(0.5)$ &$c(0.5)$ & $c(0.75)$
& $c(0.75)-c(0.5)$
\\
\hline
RQ & 22 (3) & 25 (3) & 20 (3) & 47 (9) & 52 (7) & 42 (6)
\\
EWRQ & 15 (2) & 19 (3) &16 (2) & 46 (8) & 49 (8) &40 (6)
\\
LID & 17 (2) & 18 (2) &\hphantom{0}2.9 (0.4) & 36 (5) & 42 (6) & 18 (2)
\\
Yu and Moyeed & 20 (2) & 21 (3) &13 (2) & 42 (7) &45 (6) &28 (4)
\\
Reich et al. & 20 (3) & 29 (5) & 11 (1) & \hphantom{0}4.2 (0.6)& \hphantom{0}8.6 (1.1)& \hphantom{0}3.1
(0.3)
\\
\hline
\end{tabular*}
\end{table}

From the two simulation studies, we can see that most of the time the
proposed LID method works as well as the weighted RQ, and outperforms
RQ and Yu and Moyeed's method, for estimating quantiles. LID performs
better than Reich et al.'s method in some cases and is outperformed by
Reich et al.'s method in others.
LID has a significant advantage over other methods in estimating the
difference of quantiles.
When several quantiles are of interest, including their differences,
there is a clear efficiency gain in using LID.

\subsection{Empirical studies}
\label{sec5}

In this section, we look at the June 1997 Detailed
Natality Data published by the National Center for Health
Statistics. Following the analysis in Koenker (\cite{Koe05}, page 20),
we use
65\,536 cases of recorded singleton births. We consider the following
quantile model for the birth weight data:
%
%
\begin{equation}
\label{model-birth} Q_{y_i}(\tau|x_i)=a(\tau)+b(
\tau)x_{i,1}+c(\tau)x_{i,2}+d(\tau )x_{i,3}+e(
\tau)x_{i,4},\qquad i=1,2,\ldots,n,
\end{equation}
where $x_{i,1}$ is the indicator function that indicates whether the
mother went to prenatal care for at least two times, $x_{i,2}$ is the
indicator function that indicates
whether the mother smoked or not, $x_{i,3}$ is mother's weight gain in
pounds during pregnancy, and $x_{i,4}$ is the square of mother's weight gain.
The mother's weight gain enters the model as a quadratic following the
discussion in
Koenker (\cite{Koe05}, page 23).
To make the results more comparable, we consider a slight modification
of model (\ref{model-birth}):
%
%
\begin{equation}
Q_{y_i}(\tau|x_i)=a(\tau)+b(\tau)x_{i,1}+c(
\tau)x_{i,2}+d^*(\tau )x^*_{i,3}+e^*(\tau)x^*_{i,4},
\qquad i=1,2,\ldots,n,
\end{equation}
where $x^*_{i,3}$ denotes the standardized mother's weight gain during
pregnancy and $x^*_{i,4}$ denotes the standardized square of mother's
weight gain.
We compared the results from RQ and LID (with $m=39$) for the full data
set. Here we focus on the $0.1$, $0.25$, and $0.5$ quantiles. The
results are in Table~\ref{t631}.
From the results, we can see that the estimates from both methods are
very close. The standard error from LID seems to be smaller than that
from RQ.

\begin{table}[t]
\tablewidth=\textwidth
\tabcolsep=0pt
\caption{Estimates of the parameters and their standard errors (in
parentheses) for the birth weight data}\label{t631}
{\fontsize{8}{10}\selectfont{\begin{tabular*}{\textwidth}{@{\extracolsep{\fill}}ld{2.4}d{2.4}d{2.4}d{2.4}d{2.4}d{2.4}d{2.4}d{2.4}d{2.4}d{2.4}d{2.4}d{2.4}@{\hspace*{-2pt}}}
\hline
Methods & \multicolumn{1}{l}{$b(0.1)$} & \multicolumn{1}{l}{$c(0.1)$}
& \multicolumn{1}{l}{$d^*(0.1)$} & \multicolumn{1}{l}{$e^*(0.1)$} & \multicolumn{1}{l}{$b(0.25)$} &
\multicolumn{1}{l}{$c(0.25)$} & \multicolumn{1}{l}{$d^*(0.25)$} &
\multicolumn{1}{l}{$e^*(0.25)$} & \multicolumn{1}{l}{$b(0.5)$} & \multicolumn{1}{l}{$c(0.5)$} &
\multicolumn{1}{l}{$d^*(0.5)$} & \multicolumn{1}{l@{}}{$e^*(0.5)$}
\\
\hline
RQ & -0.030 & -0.22 & 0.37 & -0.21 & -0.049 & -0.22 & 0.19 &
-0.075 & -0.061 & -0.22 & 0.127 & -0.020
\\
&(0.009)& (0.01) & (0.02) & (0.02) & (0.008) & (0.008) & (0.011) &
(0.012) &(0.006) & (0.007) & (0.008) & (0.008)
\\
LID & -0.045 & -0.22 & 0.36 & -0.22 & -0.052 & -0.23 & 0.20 &
-0.081 &-0.061 & -0.23 & 0.131 & -0.026
\\
&(0.007) & (0.003) & (0.002) & (0.003) & (0.001) & (0.002) & (0.008) &
(0.007) & (0.003) & (0.003) & (0.002) & (0.002)
\\
\hline
\end{tabular*}}}%
\end{table}

To see how good the estimates are, we compared the estimated
conditional quantile with the local quantile estimated
nonparametrically. We considered two subsets of the full data.
For the first subset of the data, we selected $x_{i,1}=1$, $x_{i,2}=1$,
and $24.5<x_{i,3}<25.5$, within which range there are 96 observations.
For the second subset of the data, we selected $x_{i,1}=1$,
$x_{i,2}=0$, and $44.5<x_{i,3}<45.5$, within which range there are 1318
observations. Then we calculated the quantile of $y_i$ in each subset
of the data as the local quantile, and compared it with the predicted
quantiles from RQ and LID. The results are presented in Table~\ref
{t634}. From the results, we can see that all the estimated quantiles
are very close to the local quantile estimates.

\begin{table}[b]
\tablewidth=\textwidth
\tabcolsep=0pt
\caption{Estimates of the local quantile}\label{t634}
\begin{tabular*}{\textwidth}{@{\extracolsep{4in minus 4in}}lllllll@{}}
\hline
& \multicolumn{3}{c}{$x_{i,1}=1$, $x_{i,2}=1$, and
$x_{i,3}=25$} & \multicolumn{3}{c@{}}{$x_{i,1}=1$, $x_{i,2}=0$, and
$x_{i,3}=45$}
\\[-5pt]
& \multicolumn{3}{c}{\hrulefill} & \multicolumn{3}{c@{}}{\hrulefill}
\\
Quantile & Local quantile & RQ & LID & Local quantile & RQ & LID
\\
\hline
0.1 & 2.54 & 2.44 & 2.43 & 2.89 & 2.90 & 2.88
\\
0.25 & 2.81 & 2.76 & 2.75& 3.18 & 3.17 & 3.17
\\
0.5 & 3.02 & 3.07 & 3.07 & 3.54 & 3.47 & 3.46
\\
\hline
\end{tabular*}
\end{table}

Another way to check the model fitness is to build the model by leaving
out a portion of the data, and then evaluate the model performance on
the out-of-bag portion of the data. Here we compared the out-of-bag
quantile coverage (the percentage of the testing data that fall below
the $\tau$th quantile line) by randomly selecting $10\%$ of the data
as the out-of-bag testing data and using the rest as the training data.
The results based on a random splitting are
summarized in Table~\ref{t102}. We can see that both RQ and LID have
coverages similar to the nominal values.

\begin{table}[t]
\tablewidth=\textwidth
\tabcolsep=0pt
\caption{Out-of-bag quantile coverage}\label{t102}
\begin{tabular*}{\textwidth}{@{\extracolsep{\fill}}llllll@{}}
\hline
Methods & $\tau= 0.1$ & $\tau= 0.25$ & $\tau= 0.5$ & $\tau=
0.75$ & $\tau= 0.9$
\\
\hline
RQ & 0.100 & 0.251 & 0.504 & 0.749 & 0.895
\\
LID & 0.093 & 0.249 & 0.506 & 0.748 & 0.909
\\
\hline
\end{tabular*}
\end{table}

From this example we can see that the model parameter estimates,
including the quantiles, from both RQ and LID are very similar, but LID
estimates are associated with lower standard errors, which corroborates
our findings in simulation studies.

\section{Conclusion}
\label{sec6}

In this paper we proposed a Bayesian method for quantile regression
which estimates multiple quantiles simultaneously. We proved the
convergence of the proposed algorithm, i.e., the stationary
distribution of the Markov chain constructed by LID would converge to
the target distribution as the number of quantiles $m$ goes to
infinity. In the simulation studies, we found that choosing $m=15$
already gave satisfactory results.
In the comparison of the proposed LID method with other methods, LID
provides comparable results for quantile estimation, and gives much
better estimates of the difference of the quantiles than other methods
(RQ, weighted RQ, and Yu and Moyeed's method).

The LID method is computationally intensive, and it requires longer
time than other methods to obtain the results. Therefore, it is of
interest to optimize LID to reduce the computational cost.

The LID method uses $m$ quantiles to construct an approximation to the
likelihood through linear interpolation.
For large $m$, it would be useful to impose regularization to make
inference more efficient. We may assume that $\beta(\tau)$ can be
characterized by a few parameters,
so we have a low-dimensional parameter space no matter what $m$ is, and
the computation of LID would simplify.
On the other hand, this approach involves additional assumption or
approximation which
would require additional work for its theoretical justification.

\begin{appendix}

\section*{Appendix: Technical details}

\subsection{Find the bounds for the proposal distribution}
\label{bound}

This is for step 3 of the algorithm in Section~\ref{algorithm}.
For each observation $(y_i,x_i)$, $i=1,2,\ldots,n$, we
can calculate a lower bound $l_{j,l,i}$ and an upper bound $u_{j,l,i}$.
Then $ l_{j,l}=\max_{i}(l_{j,l,i})$ is
taken as the maximum of all these lower bounds and $
u_{j,l}=\min_{i}(u_{j,l,i})$ is taken as the minimum of all these
upper bounds. The formula to
calculate $l_{j,l,i}$ and $u_{j,l,i}$ is given as follows.

If $1<j<m$ and $x_{i,l}>0$, where $x_{i,l}$ denotes the $l$th
element of $x_{i}$, then
\begin{eqnarray*}
&&l_{j,l,i}=\frac{x_i^{T}\beta^{k-1}(\tau_{j-1})-\sum_{t\neq l}
x_{i,t}\beta^{k-1}_t(\tau_j)}{x_{i,l}} \quad\mbox{and}
\\
&&  u_{j,l,i}=
\frac{x_i^{T}\beta^{k-1}(\tau_{j+1})-\sum_{t\neq l}
x_{i,t}\beta^{k-1}_t(\tau_j)}{x_{i,l}}.
\end{eqnarray*}
If $1<j<m$ and $x_{i,l}<0$, then
\begin{eqnarray*}
&&l_{j,l,i}=\frac{x_i^{T}\beta^{k-1}(\tau_{j+1})-\sum_{t\neq l}
x_{i,t}\beta^{k-1}_t(\tau_j)}{x_{i,l}} \quad\mbox{and}
\\
&& u_{j,l,i}=
\frac{x_i^{T}\beta^{k-1}(\tau_{j-1})-\sum_{t\neq l}
x_{i,t}\beta^{k-1}_t(\tau_j)}{x_{i,l}}.
\end{eqnarray*}
If $j=1$ and $x_{i,l}>0$, then
\[
l_{j,l,i}=-\infty\quad\mbox{and}\quad u_{j,l,i}=
\frac
{x_i^{T}\beta
^{k-1}(\tau_{j+1})-\sum_{t\neq l}
x_{i,t}\beta^{k-1}_t(\tau_j)}{x_{i,l}}.
\]
If $j=1$ and $x_{i,l}<0$, then
\[
l_{j,l,i}=\frac{x_i^{T}\beta^{k-1}(\tau_{j+1})-\sum_{t\neq l}
x_{i,t}\beta^{k-1}_t(\tau_j)}{x_{i,l}} \quad\mbox{and}\quad u_{j,l,i}=
\infty.
\]
If $j=m$ and $x_{i,l}>0$, then
\[
l_{j,l,i}=\frac{x_i^{T}\beta^{k-1}(\tau_{j-1})-\sum_{t\neq l}
x_{i,t}\beta^{k-1}_t(\tau_j)}{x_{i,l}} \quad\mbox{and}\quad u_{j,l,i}=
\infty.
\]
If $j=m$ and $x_{i,l}<0$, then
\[
l_{j,l,i}=-\infty\quad\mbox{and}\quad u_{j,l,i}=
\frac
{x_i^{T}\beta
^{k-1}(\tau_{j-1})-\sum_{t\neq l}
x_{i,t}\beta^{k-1}_t(\tau_j)}{x_{i,l}}.
\]
If $x_{i,l}=0$, then
\[
l_{j,l,i}=-\infty\quad\mbox{and}\quad u_{j,l,i}=\infty.
\]

\subsection{Proof of Theorem \texorpdfstring{\protect\ref{theo}}{3.1}}
\label{proof1}

We will verify the detailed balance condition to
show that the stationary distribution is $\hat{p}_m(B_m|X,Y)$. Denote
the probability of moving from $B_m$ to $B_m'$ by $K(B_m\rightarrow
B_m')$ and the proposal distribution
by $q(B_m\rightarrow B_m')$. We have
\begin{eqnarray*}
&&\hat{p}_m(B_m|X,Y)K\bigl(B_m\rightarrow
B_m'\bigr)
\\
&&\quad= \hat{p}_m(B_m|X,Y)q\bigl(B_m
\rightarrow B_m'\bigr)\min \biggl(1,\frac{\pi
_m(B_m'|X)\hat{p}_m(Y|X,B_m')q(B_m'\rightarrow B_m)}{\pi_m(B_m|X)\hat
{p}_m(Y|X,B_m)q(B_m\rightarrow B_m')}
\biggr)
\\
&&\quad = \frac{\pi_m(B_m|X)\hat{p}_m(Y|X,B_m)}{\hat
{p}_m(Y|X)}q\bigl(B_m\rightarrow
B_m'\bigr)\min \biggl(1,\frac{\pi_m(B_m'|X)\hat
{p}_m(Y|X,B_m')q(B_m'\rightarrow B_m)}{\pi_m(B_m|X)\hat
{p}_m(Y|X,B_m)q(B_m\rightarrow B_m')} \biggr)
\\
&&\quad = \frac{\pi_m(B_m'|X)\hat{p}_m(Y|X,B_m')}{\hat
{p}_m(Y|X)}q\bigl(B_m'\rightarrow
B_m\bigr)\min \biggl(\frac{\pi_m(B_m|X)\hat
{p}_m(Y|X,B_m)q(B_m\rightarrow B_m')}{\pi_m(B_m'|X)\hat
{p}_m(Y|X,B_m')q(B_m'\rightarrow B_m)},1 \biggr)
\\
&&\quad = \hat{p}_m\bigl(B_m'|X,Y\bigr)K
\bigl(B_m'\rightarrow B_m\bigr).
\end{eqnarray*}
So the detailed balance condition is satisfied.

\subsection{Proof of Theorem \texorpdfstring{\protect\ref{theo1}}{3.2}}
\label{proof2}

To prove Theorem~\ref{theo1}, we need three lemmas.

\begin{lemma}
\label{lemma1}
Let $\hat{p}_{m}(y_i|\theta_{m,i})=\hat{f}_{i}(y_i|x_i,B_m)$ given in
(\ref{den}). Assume $\tau_{j+1}-\tau_j=\mathrm{O}(\frac{1}{m})$.
Then
\begin{enumerate}[(b)]
\item[(a)] $|\hat{p}_{m}(y_i|\theta_{m,i})-p(y_i|\theta
_{f_i})|=\mathrm{O}(\frac
{1}{\sqrt{m}})$
uniformly in the support of $y$ as well as uniformly in~$i$.
%
\item[(b)] $|\hat{p}_m(Y|X,B_m)-p(Y|X,B_m)|=\mathrm{O}(\frac{1}{\sqrt{m}})$
uniformly in the support of $Y$.
\end{enumerate}
\end{lemma}

\begin{pf}(a)
We will prove this proposition in two different cases.

Case 1: If $y_i$ is between two quantiles we are using, in which
case we can find two consecutive quantiles $q_{i,\tau_{j}}$ and
$q_{i,\tau_{j+1}}$ such that $y_i\in[q_{i,\tau_j},q_{i,\tau_{j+1}})$,
where $1\le j\le m-1$, then by the mechanism of linear
interpolation, we have the following equation
\begin{eqnarray*}
\hat{p}_m(y_i|\theta_{m,i}) &=&
\frac{\tau_{j+1}-\tau
_{j}}{q_{i,\tau
_{j+1}}-q_{i,\tau_j}}
\\
&=& \frac{\tau_{j+1}-\tau_{j}}{F_i^{-1}(\tau_{j+1})-F_i^{-1}(\tau
_j)}
\\
&=&\frac{\tau_{j+1}-\tau_{j}}{(F_i^{-1})'(\tau^*)(\tau_{j+1}-\tau
_{j})}
\\
&=&\frac{\tau_{j+1}-\tau_{j}}{(\fraca{1}{f_i(y_i^*)})(\tau_{j+1}-\tau
_{j})}
\\
&=&f_i\bigl(y_i^*\bigr),
\end{eqnarray*}
where $\tau^*\in[\tau_{j},\tau_{j+1}),y_i^*\in
[q_{i,\tau_j},q_{i,\tau_{j+1}})$, $F_i$ denotes the cdf of
$y_i|\theta_f$, $F_i(y_i^*)=\tau^*$, and $f_i$ denotes the pdf of
$y_i|\theta_f$.

Now we want to show that
%
%
\begin{equation}
\bigl |f_i\bigl(y_i^*\bigr)-f_i(y_i)\bigr |
\leq\sup_{y\in
[q_{i,\tau_j},q_{i,\tau_{j+1}})}f_i(y)-\inf_{y\in
[q_{i,\tau_j},q_{i,\tau_{j+1}})}f_i(y)
\leq M_2\delta,
\end{equation}
where $\delta=\sqrt{2(\tau_{j+1}-\tau_{j})/M_2}$ and $M_2$ is given in
Assumption~\ref{assump1}. If
$q_{i,\tau_{j+1}}-q_{i,\tau_j}\leq\delta$, then
$|f_i(y_i^*)-f_i(y_i)|=|f_i'(y^{\dag})(y_i^*-y_i)|\leq M_2\delta$,
where $y^{\dag}\in[q_{i,\tau_j},q_{i,\tau_{j+1}})$. Now let us
consider the case that $q_{i,\tau_{j+1}}-q_{i,\tau_j}>\delta$. We
will show that
%
%
\begin{equation}
\int_{q_{i,\tau_j}}^{q_{i,\tau_{j+1}}}f_i(y)\,\mathrm{d}y>
\tau_{j+1}-\tau_{j},
\end{equation}
if
%
%
\begin{equation}
\sup_{y\in[q_{i,\tau_j},q_{i,\tau_{j+1}})}f_i(y)-\inf_{y\in
[q_{i,\tau_j},q_{i,\tau_{j+1}})}f_i(y)>
M_2\delta.
\end{equation}
Letting
$y_{\mathrm{inf}}=\arginf_{y\in[q_{i,\tau_j},q_{i,\tau_{j+1}})}f_i(y)$,
$y_{\mathrm{sup}}=\argsup_{y\in[q_{i,\tau_j},q_{i,\tau_{j+1}})}f_i(y)$,
without loss of generality, we can assume that $y_{\mathrm{inf}}<y_{\mathrm{sup}}$. It
is obvious that $y_{\mathrm{sup}}-y_{\mathrm{inf}}>\delta$, because if
$y_{\mathrm{sup}}-y_{\mathrm{inf}}\leq\delta$, then
%
%
\begin{eqnarray}
\sup_{y\in[q_{i,\tau_j},q_{i,\tau_{j+1}})}f_i(y)-\inf_{y\in
[q_{i,\tau_j},q_{i,\tau
_{j+1}})}f_i(y)&=&f_i(y_{\mathrm{sup}})-f_i(y_{\mathrm{inf}})
\nonumber
\\[-8pt]
\\[-8pt]
&=&\bigl |f_i'
\bigl(y^{\dag
}\bigr)\bigr |(y_{\mathrm{sup}}-y_{\mathrm{inf}})\leq
M_2\delta.\nonumber
\end{eqnarray}
We can find a line with slope $M_2$ that goes through
$(y_{\mathrm{sup}},f_i(y_{\mathrm{sup}}))$. This line would be below the curve
$f_i(y)$ in $[y_{\mathrm{inf}},y_{\mathrm{sup}})$, since
$f_i(y)-f_i(y_{\mathrm{sup}})=f_i'(y^{\dag\dag})(y-y_{\mathrm{sup}})\geq M_2(y-y_{\mathrm{sup}})$
for $y<y_{\mathrm{sup}}$,
which leads to $f_i(y)\geq f_i(y_{\mathrm{sup}})+M_2(y-y_{\mathrm{sup}})$.

%
\renewcommand{\thefigure}{1}
\begin{figure}[b]

\includegraphics{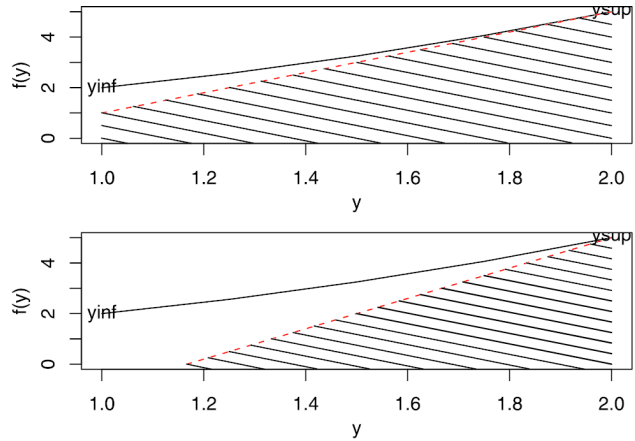}
\caption{Illustration of the two possible cases
of the area $S$:
trapezoid and triangle. The solid curve stands for $f(y)$. The dotted
line stands for
the line with slope $M_2$. The shaded area is $S$.}\label{fig1}
\end{figure}

Now we can check
the area $S$ formed by the line, $y=y_{\mathrm{inf}}$, $y=y_{\mathrm{sup}}$, and
$f_i(y)=0$. Figure~\ref{fig1} shows two possible cases. The shaded
region is $S$.

If $f_i(y_{\mathrm{sup}})-M_2(y_{\mathrm{sup}}-y_{\mathrm{inf}})\geq0$, the area is equal to
%
%
\begin{eqnarray}
\frac{[2f_i(y_{\mathrm{sup}})-M_2(y_{\mathrm{sup}}-y_{\mathrm{inf}})](y_{\mathrm{sup}}-y_{\mathrm{inf}})}{2}&\geq& \frac
{f_i(y_{\mathrm{sup}})(y_{\mathrm{sup}}-y_{\mathrm{inf}})}{2}
\nonumber
\\[-8pt]
\\[-8pt]
&>&\frac{M_2\delta^2}{2}=
\tau_{j+1}-\tau_j.\nonumber
\end{eqnarray}
If $f_i(y_{\mathrm{sup}})-M_2(y_{\mathrm{sup}}-y_{\mathrm{inf}})<0$, the area is equal to
%
%
\begin{equation}
\frac{f_i(y_{\mathrm{sup}})^2}{2M_2}>\frac{(M_2\delta)^2}{2M_2}=\tau _{j+1}-\tau_j.
\end{equation}
Therefore, in both cases, we have
%
%
\begin{equation}
\int_{q_{i,\tau_j}}^{q_{i,\tau_{j+1}}}f_i(y)\,\mathrm{d}y\geq
\int_{y_{\mathrm{inf}}}^{y_{\mathrm{sup}}}f_i(y)\,\mathrm{d}y\geq
S>\tau_{t+1}-\tau_j,
\end{equation}
which contradicts with the fact that
$\int_{q_{i,\tau_j}}^{q_{i,\tau_{j+1}}}f_i(y)\,\mathrm{d}y=\tau_{j+1}-\tau_j$.
Hence
\begin{eqnarray*}
\bigl |f_i\bigl(y_i^*\bigr)-f_i(y_i)\bigr |
&\leq& \sup_{y\in
[q_{i,\tau_j},q_{i,\tau_{j+1}})}f_i(y)-\inf_{y\in
[q_{i,\tau_j},q_{i,\tau_{j+1}})}f_i(y)
\\
&\leq& M_2\delta
= \sqrt{2M_2(\tau_{j+1}-\tau_{j})}
\\
&=& \mathrm{O}\biggl(\frac{1}{\sqrt{m}}\biggr),
\end{eqnarray*}
given that
$\tau_{j+1}-\tau_j=\mathrm{O}(\frac{1}{m})$.

Now let us consider the second case.

Case 2: If $y_i$ is a point in the tail, which means $y_i\leq
q_{i,\tau_1}$ or $y_i>q_{i,\tau_m}$, then we have
$p(y_i|\theta_{f_i})=f_i(y_i)<\frac{c}{\sqrt{m}}$ from
Assumption~\ref{assump1}. For the tail part, we can use a truncated
normal for the interpolation so that
$\hat{p}_m(y_i|\theta_{m,i})<\frac{c}{\sqrt{m}}$. Therefore, we have
$|\hat{p}_m(y_i|\theta_{m,i})-
p(y_i|\theta_{f_i})|<\frac{2c}{\sqrt{m}}=\mathrm{O}(\frac{1}{\sqrt{m}})$.

Thus for both Cases 1 and 2, we showed $|\hat{p}_m(y_i|\theta_{m,i})-
p(y_i|\theta_{f_i})|=\mathrm{O}(\frac{1}{\sqrt{m}})$.

 (b)
Let us first show
$|\hat{p}_m(y_i|\theta_{m,i})-p(y_i|x_i, B_m)|=\mathrm{O}(\frac{1}{\sqrt{m}})$.
\begin{eqnarray*}
&& \bigl |\hat{p}_m(y_i|\theta_{m,i})-p(y_i|x_i,
B_m)\bigr |
\\
&&\quad = \biggl|\int_{f_i\in\mathscr{F}_{\theta_{m,i}}} \hat{p}_m(y_i|
\theta_{m,i})\,\mathrm{d}\Pi_{\theta_{m,i}}(f_i)-\int
_{f_i\in
\mathscr
{F}_{\theta_{m,i}}} p(y_i|\theta_{f_i})\,\mathrm{d}
\Pi_{\theta_{m,i}}(f_i) \biggr|
\\
&&\quad \leq \int_{f_i\in\mathscr{F}_{\theta_{m,i}}}\bigl |\hat {p}_m(y_i|
\theta _{m,i})-p(y_i|\theta_{f_i})\bigr |\,\mathrm{d}
\Pi_{\theta_{m,i}}(f_i)
\\
&&\quad = \mathrm{O}\biggl(\frac{1}{\sqrt{m}}\biggr).
\end{eqnarray*}
Because $\hat{p}_m(Y|X,B_m)=\prod_{i=1}^n \hat{p}_m(y_i|x_i, B_m)$ and
$p(Y|X,B_m)=\prod_{i=1}^n p(y_i|x_i, B_m)$, we can show
$|\hat{p}_m(Y|X,B_m)-p(Y|X,B_m)|=\mathrm{O}(\frac{1}{\sqrt{m}})$ simply by
induction. We will show the case with $n=2$ here.
\begin{eqnarray*}
&&\bigl  |\hat{p}_m(Y|X,B_m)-p(Y|X,B_m)\bigr |
\\
&&\quad = \bigl |\hat{p}_m(y_1|X,B_m)\hat
{p}_m(y_2|X,B_m)-p(y_1|X,B_m)p(y_2|X,B_m)\bigr |
\\
&&\quad = \bigl |\hat{p}_m(y_1|X,B_m)
\hat{p}_m(y_2|X,B_m)-\hat
{p}_m(y_1|X,B_m)p(y_2|X,B_m)
\\
&&\qquad{} +\hat{p}_m(y_1|X,B_m)p(y_2|X,B_m)-p(y_1|X,B_m)p(y_2|X,B_m)\bigr |
\\
&&\quad \leq \bigl |\hat{p}_m(y_1|X,B_m)\bigl[
\hat{p}_m(y_2|X,B_m)-p(y_2|X,B_m)
\bigr]\bigr |
\\
&&\qquad{} +\bigl |\bigl[\hat{p}_m(y_1|X,B_m)-p(y_1|X,B_m)
\bigr]p(y_2|X,B_m)\bigr |
\\
&&\quad= M_1\mathrm{O}\biggl(\frac{1}{\sqrt{m}}
\biggr)+M_1\mathrm{O}\biggl(\frac{1}{\sqrt{m}}\biggr)
\\
&&\quad = \mathrm{O}\biggl(\frac{1}{\sqrt{m}}\biggr),
\end{eqnarray*}
where $M_1$ is given in Assumption~\ref{assump1}.
The proof can be easily generalized to the case with $n>2$.
\end{pf}

\begin{lemma}
\label{lemma2}
\begin{enumerate}[(b)]
\item[(a)] $E_{\pi_m}(|\hat{p}_m(Y|X,B_m)-p(Y|X,B_m)|)=\mathrm{O}(\frac
{1}{\sqrt{m}})$.
\item[(b)] $E_{\pi_m}(|\hat{p}_m(Y|X,B_m)-\hat
{p}_{m-1}(Y|X,B_{m-1})|)=\mathrm{O}(\frac{1}{\sqrt{m}})$.
\end{enumerate}
\end{lemma}

\begin{pf}
Part (a) of Lemma~\ref{lemma2} follows immediately from Lemma~\ref
{lemma1}(b). Part (b) of Lemma~\ref{lemma2} can be obtained by applying
Lemma~\ref{lemma2}(a) twice.
\end{pf}

\begin{lemma}
\label{lemma3}
$|\hat{p}_m(Y|X)-p(Y|X)|=\mathrm{O}(\frac{1}{\sqrt{m}})$.
\end{lemma}

\begin{pf}
\begin{eqnarray*}
&& \bigl |\hat{p}_m(Y|X)-p(Y|X)\bigr |
\\
&&\quad = \biggl\llvert \int\pi_m(B_m|X)\bigl[
\hat{p}_m(Y|X,B_m)-p(Y|X,B_m)\bigr]\,
\mathrm{d}B_m\biggr\rrvert
\\
&&\quad \leq \int\pi_m(B_m|X)|\hat{p}_m(Y|X,B_m)-p(Y|X,B_m)|
\,\mathrm{d}B_m
\\
&&\quad = E_{\pi_m}\bigl(\bigl |\hat{p}_m(Y|X,B_m)-p(Y|X,B_m)\bigr |
\bigr)
\\
&&\quad = \mathrm{O}\biggl(\frac{1}{\sqrt{m}}\biggr).
\end{eqnarray*}
\upqed\end{pf}

Now we are ready to prove Theorem~\ref{theo1}. We have
\begin{eqnarray*}
\hspace*{-2pt}&& \bigl \Vert \hat{p}_m(B_m|X,Y)-p(B_m|X,Y)
\bigr \Vert _{\mathrm{TV}}
\\
\hspace*{-2pt}&&\quad = \frac{1}{2}\sup_{|h|\leq1}\biggl\llvert \int
h(B_m) \biggl(\frac{\pi
_m(B_m|X)\hat{p}_m(Y|X,B_m)}{\hat{p}_m(Y|X)}-\frac{\pi
_m(B_m|X)p(Y|X,B_m)}{p(Y|X)} \biggr)\,
\mathrm{d}B_m\biggr\rrvert
\\
\hspace*{-2pt}&&\quad \leq \frac{1}{2}\int\pi_m(B_m|X)\biggl
\llvert \frac{\hat
{p}_m(Y|X,B_m)}{\hat
{p}_m(Y|X)}-\frac{p(Y|X,B_m)}{p(Y|X)}\biggr\rrvert \,\mathrm{d}B_m
\\
\hspace*{-2pt}&&\quad= \frac{1}{2}\int\pi_m(B_m|X)\biggl
\llvert \frac{\hat
{p}_m(Y|X,B_m)p(Y|X)-\hat
{p}_m(Y|X)p(Y|X,B_m)}{\hat{p}_m(Y|X)p(Y|X)}\biggr\rrvert \,\mathrm{d}B_m
\\
\hspace*{-2pt}&&\quad = \frac{1}{2}\int\pi_m(B_m|X)
\\
\hspace*{-2pt}&&\phantom{\quad = \frac{1}{2}\int}{}\times\biggl\llvert \frac{[\hat
{p}_m(Y|X,B_m)-p(Y|X,B_m)]p(Y|X)+p(Y|X,B_m)[p(Y|X)-\hat
{p}_m(Y|X)]}{\hat
{p}_m(Y|X)p(Y|X)}\biggr\rrvert \,\mathrm{d}B_m
\\
\hspace*{-2pt}&&\quad\leq \frac{1}{2}\int\pi_m(B_m|X)
\\
\hspace*{-2pt}&&\phantom{\quad = \frac{1}{2}\int}{}\times\frac{|\hat
{p}_m(Y|X,B_m)-p(Y|X,B_m)|p(Y|X)+p(Y|X,B_m)|p(Y|X)-\hat
{p}_m(Y|X)|}{\hat
{p}_m(Y|X)p(Y|X)}\,\mathrm{d}B_m
\\
\hspace*{-2pt}&&\quad = \frac{1}{2} \biggl[\frac{E_{\pi_m}(|\hat
{p}_m(Y|X,B_m)-p(Y|X,B_m)|)}{\hat{p}_m(Y|X)}+\frac{|\hat
{p}_m(Y|X)-p(Y|X)|}{\hat{p}_m(Y|X)}
\biggr].
\end{eqnarray*}
We already know from Lemma~\ref{lemma3} that $\hat
{p}_m(Y|X)\rightarrow
p(Y|X)$ as
$m\to\infty$, so for any $e^*\in(0, p(Y|X))$,
there exists an $m^*$ such
that $|\hat{p}_m(Y|X)-p(Y|X)|<e^*$ for $m\geq m^*$. We can see that
\[
LB=\min\bigl(\hat{p}_{m_0}(Y|X),\hat{p}_{m_0+1}(Y|X),\ldots,\hat
{p}_{m^*-1}(Y|X),p(Y|X)-e^*\bigr)
\]
is a lower bound for $\hat{p}_m(Y|X)$, where $m_0$ is the minimum
number of quantiles we use.
Therefore, $\Vert \hat{p}_m(B_m|X,Y)-
p(B_m|X,Y)\Vert _{\mathrm{TV}}\leq\frac{1}{2LB}[E_{\pi_m}(|\hat
{p}_m(Y|X,B_m)-p(Y|X,B_m)|)+|\hat{p}_m(Y|X)-p(Y|X)|]=\mathrm{O}(\frac{1}{\sqrt {m}})\rightarrow0$
as $m\rightarrow\infty$ (because of Lemmas~\ref{lemma2} and~\ref{lemma3}).
\end{appendix}

\section*{Acknowledgements}

The research of Yuguo Chen was supported in part by NSF Grant
DMS-1106796. The research
of Xuming He was supported in part by NSF Grants DMS-1237234 and
DMS-1307566, and National Natural Science Foundation of China Grant 11129101.


%

\printhistory

\end{document}